\documentclass[11pt, a4paper]{article}

\usepackage[utf8]{inputenc}
\usepackage[T1]{fontenc}
\usepackage[french, english]{babel} 
\usepackage{amsmath, amssymb}
\usepackage{graphicx}
\usepackage{hyperref}
\usepackage{geometry}
\usepackage{xcolor, colortbl}
\usepackage{booktabs}
\usepackage{caption}
\usepackage{subcaption}
\usepackage{mathrsfs}
\usepackage{authblk}
\usepackage{tikz}
\usetikzlibrary{positioning, arrows.meta, shapes.geometric, calc, backgrounds, fit}

\newcommand{\cD}[0]{\mathcal{D}}

\newcommand{\DFD}{\cD_{\mathrm{FD}}}
\newcommand{\DMC}{\cD_{\mathrm{MC}}}

\title{Data-efficient extraction of optical properties from 3D Monte Carlo TPSFs using Bi-LSTM transfer learning}
\author[1,2]{Joubine Aghili}
\author[2]{R\'{e}mi Imbach}
\author[3,4]{A. Pallarès}
\author[3]{P. Schmitt}
\author[3]{W. Uhring}

\affil[1]{\small IRMA, Université de Strasbourg, CNRS UMR 7501, 7 rue René Descartes, 67084
Strasbourg, France}
\affil[2]{\small INRIA Nancy-Grand Est, MACARON Project, Strasbourg, France}
\affil[3]{\small Université de Strasbourg, CNRS, ENGEES, ICube UMR 7357, Strasbourg F-67000, France}
\affil[4]{\small Université de Haute Alsace, Mulhouse F-68100, France}

\date{\today}

\begin{document}

\maketitle

\begin{abstract}
  Time-Resolved Spectroscopy (TRS) is a powerful modality for non-invasive characterization of turbid media. However, extracting optical properties—absorption $\mu_a$ and reduced scattering $\mu_s'$ — from 3D stochastic measurements remains computationally expensive for real-time applications.
  In this paper, we propose a data-efficient, physics-informed transfer learning strategy using a Bidirectional Long Short-Term Memory (Bi-LSTM) network.
  By leveraging a fast deterministic solver to establish a physical prior before fine-tuning on a restricted set of 3D Monte Carlo simulations, our model successfully bridges the analytical-to-stochastic domain gap.
  The proposed method eliminates the systematic bias of analytical models while maintaining a competitive error with near-instantaneous inference time.
\end{abstract}

\vspace{1em}

\section{Introduction}
\label{sec:intro}

Accurately determining the optical properties of turbid media is a fundamental challenge with broad implications across diverse scientific fields. Applications range from biomedical imaging and tissue diagnostics \cite{jacques2013}, which have driven the development of advanced non-invasive measurement devices, to atmospheric and climate sciences \cite{fussen2025,chandrasekhar1944}. The primary parameters characterizing such highly scattering environments are the absorption coefficient ($\mu_a$) and the reduced scattering coefficient ($\mu_s'$). Among the established experimental techniques to extract these parameters, Time-Resolved Spectroscopy (TRS) stands out as one of the most informative. This approach consists of injecting an ultrashort light pulse (typically in the 100 ps range) into the medium and analyzing the temporal evolution of the exiting light. The resulting Temporal Point Spread Function (TPSF) encapsulates rich, high-dimensional information regarding the photon dynamics within the medium, as the distribution of the photons' time-of-flight often spreads over several nanoseconds.

In practice, simulating photon dynamics relies on two major classes of methods: deterministic and stochastic solvers. Deterministic approaches solve the unsteady Radiative Transfer Equation (RTE) using numerical methods like finite differences \cite{turek1993} or finite elements \cite{richling2001}. While offering rigorous error estimates and the ability to handle complex geometries, achieving high temporal resolution is computationally expensive. Conversely, stochastic Monte Carlo (MC) simulations \cite{korkin2022} track individual photon trajectories, providing a \textit{gold standard} for 3D light transport at the cost of significant computation time, making them impractical for real-time inversion.

Machine learning has emerged as a promising direction to bypass these limitations. Recent works have proposed surrogate models to accelerate the forward phase or solve the inverse problem directly \cite{garcia-cuesta2009, mishra2021}. For time-resolved data, architectures like Long Short-Term Memory (LSTM) networks have been used to invert TPSF signals \cite{li2010}. Notably, Zhang et al. (2024) demonstrated that employing complete temporal profiles yields superior accuracy compared to reduced features like mean time-of-flight \cite{zhang2024}. 

Despite these advancements, deep learning typically requires massive datasets—often exceeding 100,000 simulations—to achieve high accuracy \cite{fussen2025}. Generating sufficiently many Monte Carlo simulations is prohibitive. Furthermore, models trained exclusively on fast deterministic priors suffer from a severe \textit{domain gap} when faced with realistic 3D stochastic data, leading to systematic biases. Bridging this accuracy-speed gap without relying on exorbitant datasets remains a major open challenge.

To further understand this disparity, Principal Component Analysis (PCA) offers a compelling perspective. While deterministic models often exhibit low intrinsic dimensionality (explaining 99\% variance with few modes), 3D stochastic simulations reveal a much more complex manifold, requiring over a hundred modes. This structural gap, possibly rooted in 3D boundary losses and early sub-diffusive regimes, motivates the use of high-capacity non-linear architectures like Bi-LSTMs \cite{graves2005}.

While transfer learning has become a standard and highly successful procedure in computer vision, recent extensive studies have demonstrated that transferring deep neural network weights is also exceptionally effective for time series classification tasks, significantly improving model generalization when adapting to new target datasets \cite{fawaz2018}. Inspired by these advances, we formulate the inversion of time-resolved optical signals (TPSFs) as a sequential time series problem and propose a physics-informed transfer learning strategy.

In this paper, we first train a Bi-LSTM on a large, inexpensive deterministic dataset to acquire a physical prior, then fine-tune it on a restricted set of 3,400 MC simulations. This approach successfully bridges the domain gap while reducing data requirements. 

The remainder of this paper is organized as follows: Section \ref{sec:generation} details the forward modeling and dataset analysis. Section \ref{sec:dl} presents the deep learning framework. Section \ref{sec:results} demonstrates the results and the bridging of the domain gap, followed by conclusions in Section \ref{sec:conclusion}.

\section{Forward Modeling and Dataset Generation}
\label{sec:generation}

\subsection{Measurement Geometry and Physics}
\label{sec:geometry_physics}

\begin{figure}
  \centering
  \includegraphics[width=0.48\textwidth]{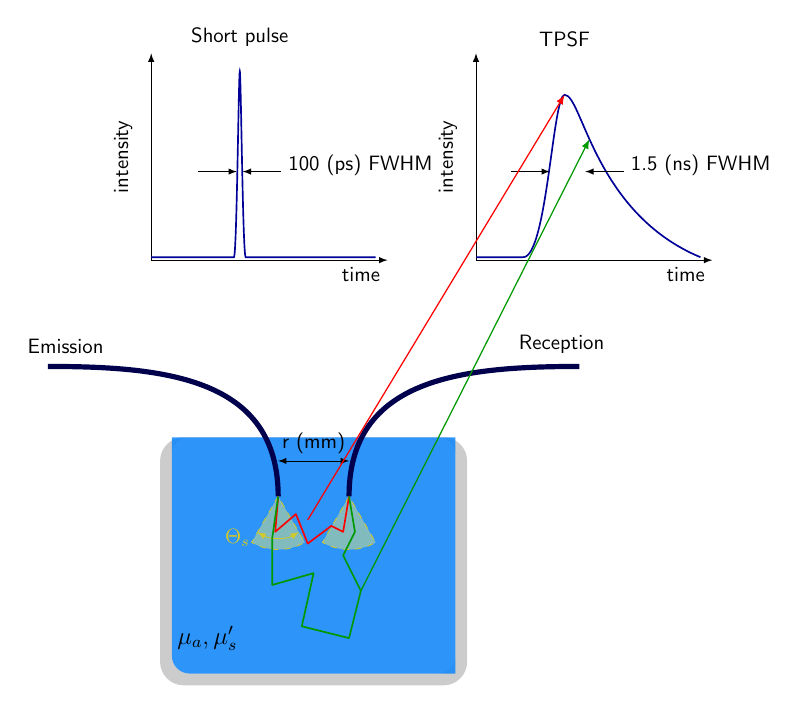}
  \includegraphics[width=0.48\textwidth]{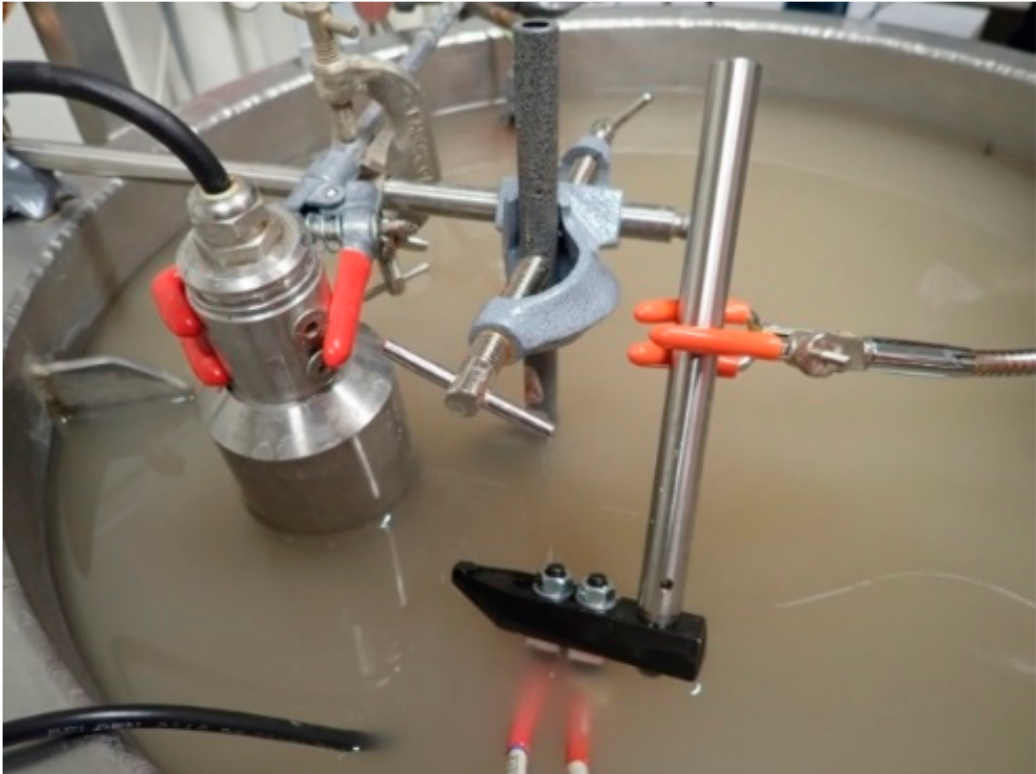}
  \caption{(Left) Schematic representation of the TROT measurement device; (Right) a photography of the device during the experiment \cite{pallares2021}}
  \label{fig:trot}
\end{figure}

The experimental foundation of this work is based on the Time-Resolved Optical Turbidity (TROT) technique \cite{pallares2021}, which aims to characterize turbid media by measuring the reflectance of light in the time domain. The core principle involves injecting an ultrashort optical pulse, with a duration of approximately 100 ps, from a laser source into the medium via an excitation fiber. A separate detection fiber, positioned at a precise and constant source-detector separation distance $\rho$, collects the scattered photons that have migrated through the medium. The exiting optical signal is captured by a high-sensitivity photo-detector and processed by a Time-Correlated Single-Photon Counting (TCSPC) system to reconstruct the measured Temporal Point Spread Function (TPSF). A simplified schematic of this setup is illustrated in Figure \ref{fig:trot}.

The propagation of light within such scattering environments is physically governed by the Radiative Transfer Equation (RTE). For a computational domain $\Omega \subset \mathbb{R}^2$ (or $\mathbb{R}^3$), the radiance $u(t, \mathbf{x}, \mathbf{s})$ at time $t$, position $\mathbf{x}$, and direction $\mathbf{s}$ is described by the following integro-differential equation:
\begin{equation}
  \label{eq:RTE_compact}
  \frac{1}{v}\partial_t u + (\mathbf{s}\cdot\nabla) u + (\mu_a + \mu_s) u = \mu_s \int_{\Theta} u(t, \mathbf{x}, \mathbf{s}') \phi_g(\mathbf{s}, \mathbf{s}') d\mathbf{s}' + S(t, \mathbf{x}, \mathbf{s})
\end{equation}
where $v$ is the speed of light in the medium, $\mu_a$ is the absorption coefficient, and $\mu_s$ is the scattering coefficient. The phase function $\phi_g$, often modeled by the Henyey--Greenstein function, describes the probability of a photon being scattered from direction $\mathbf{s}'$ to $\mathbf{s}$, and is characterized by the anisotropy factor $g$. In the study of highly scattering turbid media, it is standard practice to work directly with the reduced scattering coefficient $\mu_s' = \mu_s(1-g)$, as it effectively combines the scattering density and anisotropy into a single isotropic parameter, which is the primary variable of interest for characterizing turbidity. Thus, the inverse problem focuses on the simultaneous extraction of the pair $(\mu_a, \mu_s')$ from the high-dimensional information encoded within the TPSF.

\subsection{Deterministic Solver (Source Domain)}
\label{sec:deterministic_solver}

To establish a baseline physical prior for our neural network, we employ a deterministic solver based on the 2D Finite Difference Discrete Ordinate Method (FD-DOM). The Discrete Ordinate Method simplifies the angular dependence of the RTE by discretizing the orientation variable $\mathbf{s}$ into $M$ fixed angles on the unit circle, thereby transforming the integro-differential equation into a system of $M$ coupled partial differential equations. 

These equations are subsequently solved on a spatial Cartesian grid $(x_i, y_j)$. The spatial advection terms are approximated using a standard explicit upwind finite difference scheme, which ensures numerical stability while capturing the directional flow of photons. The time evolution is handled via explicit time-stepping ($\Delta t$), and zero incoming flux boundary conditions ($\mathbf{s} \cdot \mathbf{n} < 0$, where $\mathbf{n}$ is the outward normal vector) are enforced at the domain edges to properly simulate the physical boundaries. The source term is modeled to replicate the 100 ps short pulse injected by the experimental device.

While this 2D deterministic formulation inherently simplifies the complex 3D reality—most notably by neglecting out-of-plane scattering losses—it offers a paramount advantage: extreme computational efficiency. Because the solutions preserve the regularity of the forcing source and are free of stochastic noise, this solver provides a perfectly smooth mathematical baseline. 

We leverage this highly optimized solver to rapidly generate a large-scale, noise-free \textit{source dataset}, denoted $\DFD$. For this study, $\DFD$ comprises exactly 7,441 TPSF curves generated across a wide, uniformly sampled range of absorption and reduced scattering coefficients. This dataset serves as the foundational training ground, allowing the deep learning architecture to learn the fundamental non-linear temporal dynamics of light propagation at a fraction of the computational cost required by high-fidelity stochastic methods.

\subsection{Stochastic Solver (Target Domain)}
\label{sec:stochastic_solver}

To capture the true complexity of light transport, including 3D geometries and stochastic shot noise, we developed a high-fidelity, GPU-accelerated Monte Carlo (MC) code implemented in Python and CUDA. This stochastic solver acts as our ground truth and generates the \textit{target dataset}, denoted $\DMC$. 

Unlike the 2D deterministic approximation, the MC simulations rigorously model photon transport in a three-dimensional homogeneous turbid medium. The physical domain is defined as a cylinder with a radius of $R = 60$ mm and a height of $H = 60$ mm, with a refractive index of $n = 1.33$. The illumination source is modeled as an optical fiber (radius $r_f = 1.0$ mm) emitting a temporal Gaussian pulse ($\Delta t_\mathrm{FWHM} = 100$ ps) into a cone with a half-angle of $\theta = 30^\circ$. A detection fiber is placed at a source-detector separation distance of $d = 10$ mm to perfectly replicate the experimental TROT setup.

The algorithm launches $N = 10^7$ photon packets per simulation. To optimize computational efficiency, a weight-based variance reduction technique is applied: at each scattering event, the photon's weight is attenuated by a factor of $e^{-\mu_a \ell}$, where $\ell$ is the path length traveled between scattering events. To avoid tracking photons with negligible contributions, a Russian roulette scheme is triggered when a photon's weight drops below $w_\mathrm{min} = 10^{-3}$; the photon survives with a probability of $p = 0.1$ (its weight being multiplied by $1/p$ to ensure energy conservation) or is otherwise annihilated. 

The temporal response is recorded over a time window of $T = 1$ ns, discretized into $N_\mathrm{bins} = 200$ intervals, yielding the final stochastic TPSF. While this solver inherently accounts for realistic 3D boundary losses and the highly directional sub-diffusive regime of early photons (the theoretical ballistic peak appearing at approximately $0.244$ ns), its high computational burden limits the size of the target dataset. 
Consequently, the generated $\DMC$ dataset consists of exactly 3,700 TPSF profiles, serving as the restricted, highly realistic dataset for the fine-tuning phase of our transfer learning architecture.
\begin{figure}[htbp]
  \centering
  \includegraphics[width=1.0\textwidth]{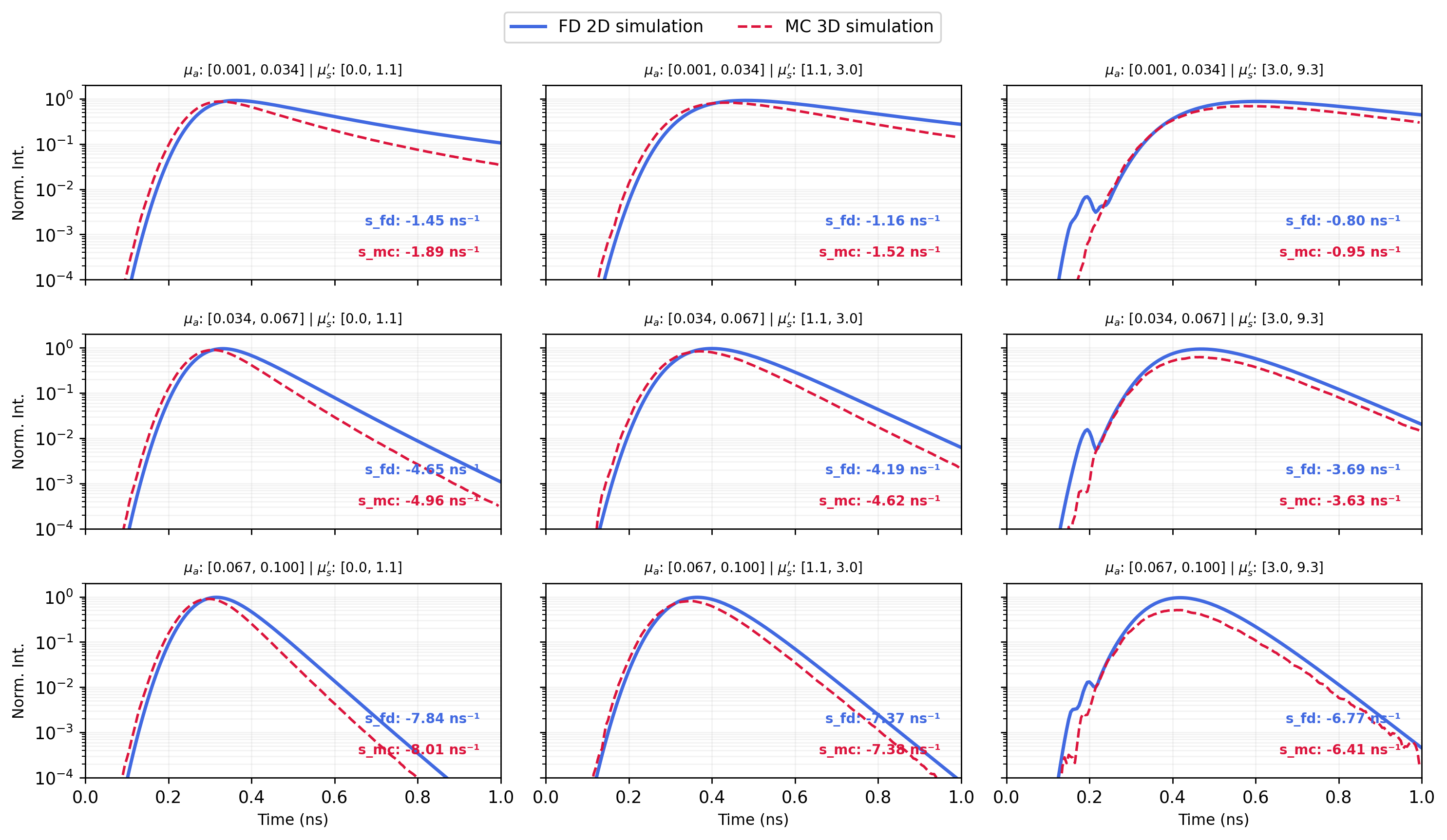}
  \caption{Comparison of the temporal decay profiles between the deterministic ($\DFD$) and stochastic ($\DMC$) datasets. The Temporal Point Spread Functions (TPSFs) are visualized on a logarithmic scale and colored by their asymptotic decay slope $m$, evaluated in the late-time regime ($t > 0.75$ ns). While the analytical FD model yields perfectly smooth exponential decays, the 3D Monte Carlo data exhibits pronounced stochastic shot noise in the tail due to the finite number of detected photon packets.}
  \label{fig:datasets.comparison}
\end{figure}
To visually underscore the fundamental physical differences between the analytical prior and the stochastic reality, Figure \ref{fig:datasets.comparison} presents a side-by-side comparison of the normalized TPSF signals from both datasets ($\DFD$ and $\DMC$).
The signals are plotted on a logarithmic scale to highlight the late-time exponential decay, with curves colored according to their asymptotic slope $m$ evaluated for $t > 0.75$ ns.
As expected from the diffusion theory, the deterministic solver produces perfectly smooth, linear tails.
In contrast, the Monte Carlo dataset intrinsically suffers from severe stochastic shot noise in this late-time regime.
As the time-of-flight increases, the probability of detecting a photon packet drops exponentially, leading to highly noisy and erratic tails.
This visual comparison explicitly illustrates the complexity of realistic optical inversion: a neural network trained solely on smooth analytical data will inevitably struggle to robustly extract optical properties—particularly the absorption coefficient, which is heavily reliant on the late-time slope—when confronted with the noisy, sub-diffusive reality of 3D stochastic transport.

\subsection{Statistical Analysis and the Domain Gap}
\label{sec:gap_analysis}

\begin{figure}[htbp]
  \centering
  \includegraphics[width=0.48\textwidth]{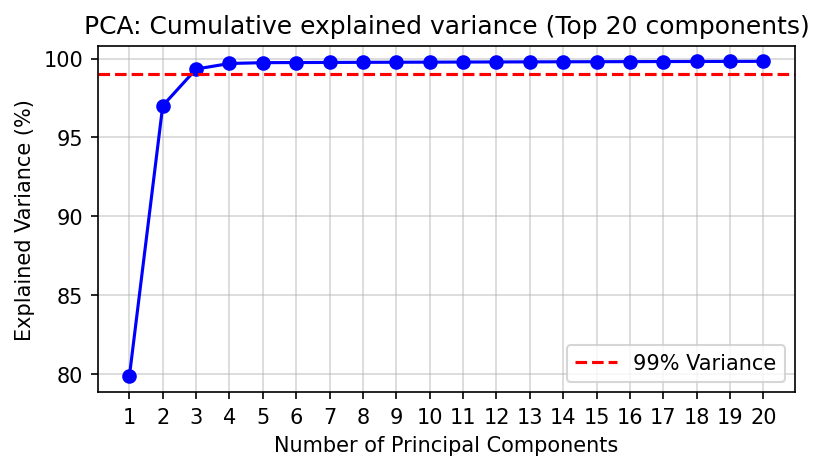}   
  \includegraphics[width=0.48\textwidth]{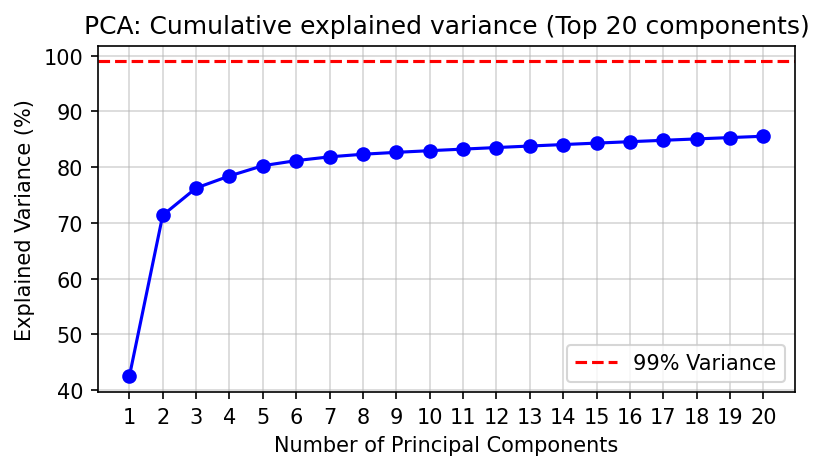} 
  \caption{Cumulative explained variance by Principal Components for both the Finite Difference (left) and Monte Carlo (right) TPSF datasets.}
  \label{fig:pca_variance}
\end{figure}

To characterize the intrinsic complexity of the generated datasets and rigorously quantify the domain gap between the deterministic approximation and the stochastic reality, a Principal Component Analysis (PCA) was performed on the log-transformed Temporal Point Spread Functions (TPSF). As illustrated in Figure \ref{fig:pca_variance}, the topological structure of the Finite Difference (FD) dataset is extremely simple: merely 3 principal components are sufficient to explain over 99\% of the signal's variance. This low dimensionality reflects the idealized, smooth nature of the diffusion equation. Conversely, the 3D Monte Carlo dataset exhibits a drastically higher intrinsic dimensionality, requiring over 130 principal components to capture 99\% of the variance. This explosion in dimensionality is not solely due to stochastic shot noise; it is deeply rooted in physical phenomena ignored by the analytical model, namely 3D boundary photon losses and the highly non-linear dynamics of the early sub-diffusive regime. This structural contrast mathematically demonstrates that simple linear inversions are obsolete for realistic photon dynamics, thereby justifying the use of high-capacity, non-linear deep learning architectures.

\begin{figure}[htbp]
  \centering
  \includegraphics[width=0.48\textwidth]{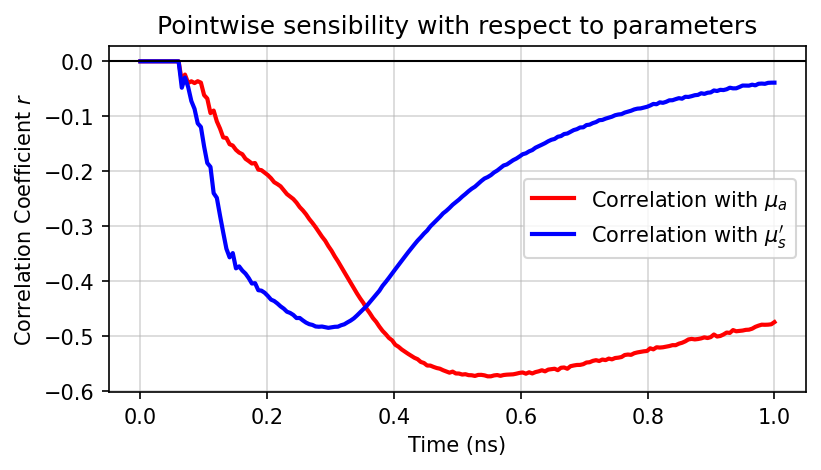}
  \includegraphics[width=0.48\textwidth]{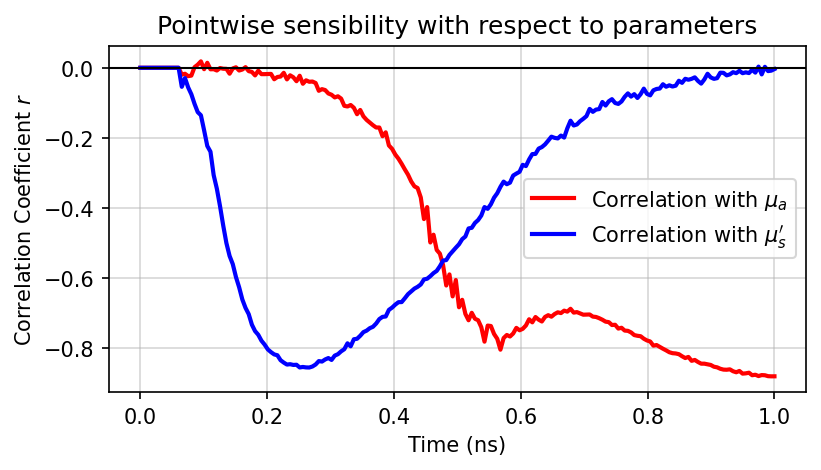}
  \caption{Temporal Pearson correlation between TPSF time series and the optical parameters ($\mu_a$, $\mu_s'$) for FD (left) and MC (right) datasets.}
  \label{fig:temporal_correlation}
\end{figure}

Furthermore, a temporal Pearson correlation analysis between the individual TPSF time bins and the target optical properties reveals critical differences in photon transport dynamics (see Figure \ref{fig:temporal_correlation}). In both datasets, the absorption coefficient $\mu_a$ exhibits a strong negative correlation on the late-time tail of the pulse, aligning with the Beer-Lambert law which exponentially attenuates photons with long path lengths. However, the influence of the reduced scattering coefficient $\mu_s'$ exposes the mechanistic limitations of the analytical model at early times. In the FD dataset, the correlation with $\mu_s'$ peaks immediately and sharply at the signal's onset, reflecting the mathematical assumption of an instantly isotropic point source. In stark contrast, the MC dataset displays a delayed and broadened correlation peak for $\mu_s'$.
This delay captures the sub-diffusive regime, where early ballistic and ``snake'' photons retain a strong directional inertia before scattering events force them into a fully diffuse state.

These fundamental temporal discrepancies dictate that a neural network trained solely on theoretical FD data will systematically misinterpret the early-time features of actual stochastic signals. This inherent \textit{domain gap} directly motivates our physics-informed transfer learning strategy, which leverages the FD data for rapid prior acquisition and the MC data for precise fine-tuning.

\section{Deep Learning Framework}
\label{sec:dl}

To overcome the inherent high-dimensionality and non-linearity of the stochastic photon transport described in Section \ref{sec:generation}, we propose a deep learning framework designed to process the sequential nature of the TPSF. The approach relies on a Bidirectional Long Short-Term Memory (Bi-LSTM) network, optimized through a two-stage physics-informed transfer learning strategy.

\subsection{Network Architecture: Dual-Head Bi-LSTM}
\label{sec:architecture}

The Temporal Point Spread Function is fundamentally a time-series signal where different temporal regions encode different physical properties. As demonstrated in our correlation analysis (Section \ref{sec:gap_analysis}), the early-time regime is heavily influenced by scattering ($\mu_s'$), whereas the late-time exponential decay is primarily driven by absorption ($\mu_a$).

To fully exploit this temporal separation, we employ a Bidirectional LSTM \cite{graves2005}. The forward pass captures the causal propagation of light from the ballistic onset to the diffuse tail, while the backward pass effectively reads the exponential decay in reverse, providing a strong inductive bias for extracting the absorption coefficient.

\begin{figure}[htbp]
  \centering
  \includegraphics[width=\textwidth]{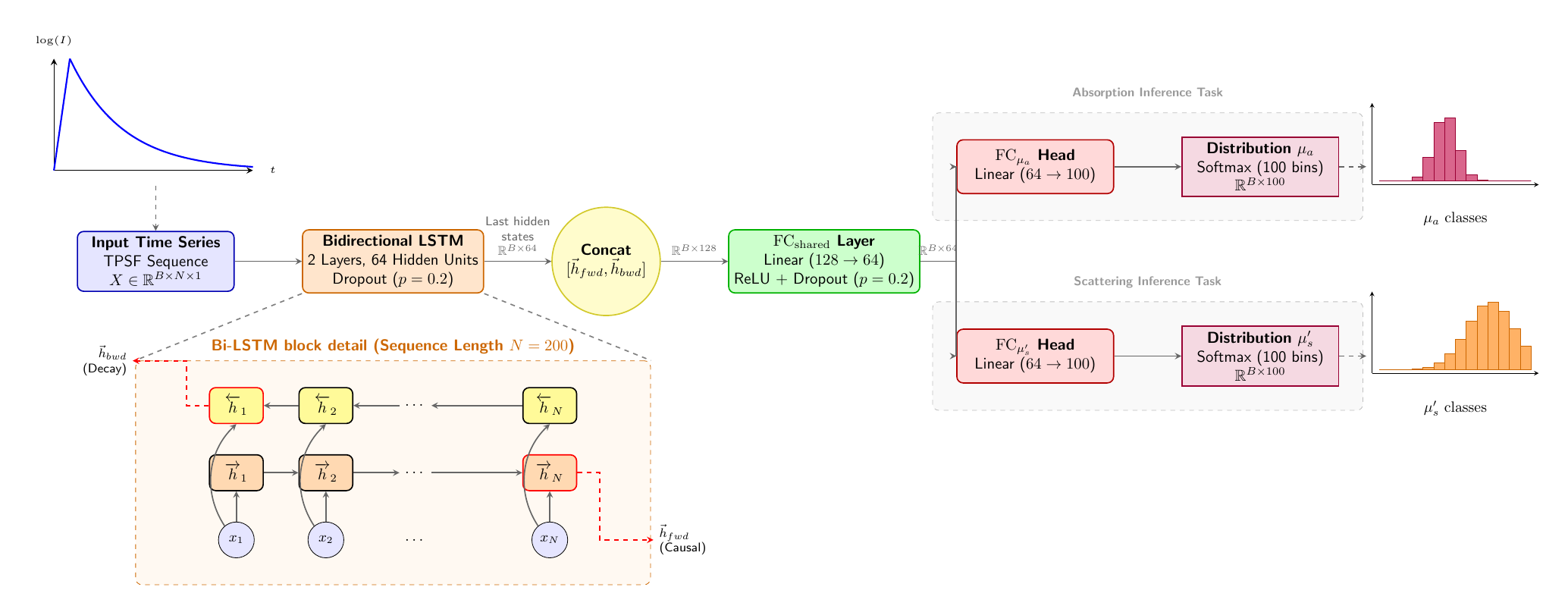}
  \caption{Architecture of the proposed Dual-Head Bi-LSTM. The normalized TPSF is processed bidirectionally. The final hidden state is routed to two independent fully connected heads configured for classification, predicting the probability distributions of $\mu_a$ and $\mu_s'$ simultaneously.}
  \label{fig:architecture}
\end{figure}

To mitigate the well-documented issue of physical crosstalk—where the network struggles to untangle the coupled effects of absorption and scattering—we branch the final hidden state into a ``Dual-Head'' configuration. Crucially, instead of performing a standard scalar regression, we reframe the inversion as a \textit{dual-classification task}. The continuous physical parameter spaces of $\mu_a$ and $\mu_s'$ are discretized into a predefined number of classes (bins).
Two independent Fully Connected (FC) blocks are used to output a probability distribution over these bins via a Softmax activation, see Figure \ref{fig:architecture}.
This explicit separation prevents the gradients of one parameter from interfering with the other during backpropagation, while simultaneously providing a direct visualization of the predictive uncertainty.

In view of Figure \ref{fig:temporal_correlation}, extensive preliminary testings with more complex mechanisms, such as temporal Dual-Attention, revealed a tendency to overfit the stochastic shot noise inherent to the small Monte Carlo dataset.
The standard Bi-LSTM architecture, shown in Figure \ref{fig:architecture}, proved to be the most robust feature extractor against stochastic noise.

\subsection{Physics-Informed Transfer Learning Methodology}
\label{sec:transfer_learning}

Training a deep recurrent neural network directly on the complex, noisy target dataset $\DMC$ (3,700 samples) would inevitably lead to severe overfitting. To bridge the domain gap, we deploy a two-step training methodology.

\textbf{Step 1: Physical Prior Acquisition (Pre-training).} \\
In the first phase, the network is trained from scratch entirely on the large, noise-free analytical dataset $\DFD$. Using a standard 80/20 split, the model learns the theoretical prior from approximately 5,952 deterministic samples. The objective is for the Bi-LSTM network (see \ref{fig:architecture}) to learn the fundamental, idealized relationships between the TPSF shape and the optical parameters. The network is trained using an Adam optimizer and a composite cross-entropy loss function. By reframing the problem as classification, the network learns to confidently predict the correct parameter bin. At the end of this phase, the model achieves high accuracy on theoretical data, acting as a robust mathematical baseline.

\textbf{Step 2: Pre-processing and Bridging the Domain Gap (Fine-Tuning).} \\
In the second phase, the pre-trained weights are loaded, and the model is fine-tuned using the stochastic target dataset $\DMC$. To rigorously evaluate the model and prevent data leakage, the 3,700 stochastic samples are split into a training set (2,960 samples) and a hold-out test set (740 samples). 

To maximize data efficiency and network robustness, a data augmentation strategy is employed: proportional Gaussian noise is added to the raw MC signals to artificially double the training set size, yielding 5,920 augmented training samples. Crucially, before feeding these augmented stochastic signals to the network, a Savitzky-Golay filter \cite{savitzky1964} is systematically applied to the inputs. This specific filtering operation smooths the high-frequency stochastic shot noise while perfectly preserving the crucial low-frequency physical dynamics, such as the exact position of the ballistic peak and the slope of the exponential decay.

To preserve the physical prior while allowing adaptation to the 3D reality, we apply a differential learning rate strategy using the Adam optimizer. The learning rate of the core Bi-LSTM layers ($\texttt{lr\_lstm}$) is strictly decoupled from the rest of the network. These recurrent layers are either trained with a heavily reduced learning rate or completely frozen ($\texttt{freeze\_lstm}$), effectively acting as a robust feature extractor. Conversely, the fully connected classification blocks retain a higher learning rate ($\texttt{lr\_fc}$). This targeted flexibility forces the network to recalibrate its probability distributions, effectively shifting its temporal focus to account for the delayed sub-diffusive regime present in the MC data. Furthermore, a \texttt{ReduceLROnPlateau} scheduler is employed to dynamically halve the learning rates if the validation loss stagnates, ensuring precise convergence.

\section{Results and Discussion}
\label{sec:results}

To rigorously evaluate the proposed methodology, we utilize several standard regression metrics: the Mean Relative Error (MRE), the Mean Error (Bias), the Standard Deviation (Std) of the relative errors, and the Success Rate, defined as the percentage of predictions falling within a 10\% error margin relative to the ground truth.

\subsection{Baseline Performance on the Source Domain}
\label{sec:res_baseline}

Initially, the Bi-LSTM network was trained and evaluated strictly within the deterministic source domain ($\DFD$). Evaluated on a hold-out test set of 1,489 analytical TPSF curves, the model demonstrated exceptional reconstruction capabilities.

\begin{figure}[htbp]
  \centering
  \includegraphics[width=1.0\textwidth]{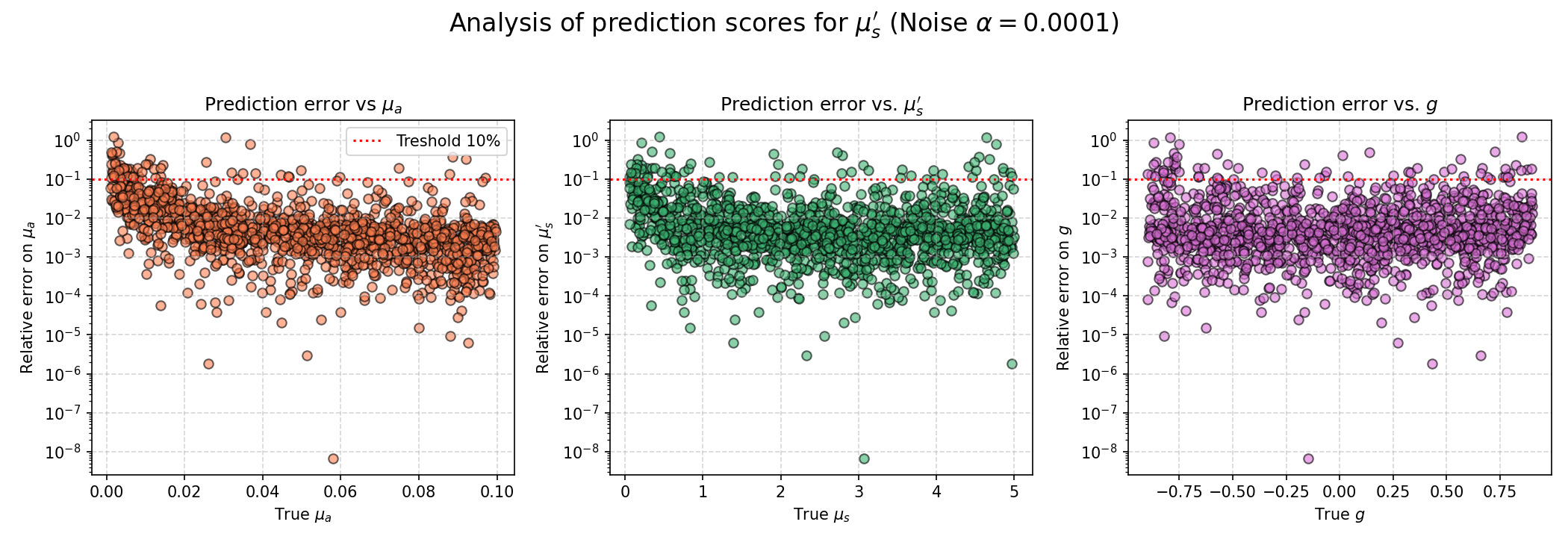}
  \caption{Baseline model predictions evaluated on the deterministic $\DFD$ hold-out test set. The strong alignment with the ground truth (solid line) demonstrates the network's capacity to perfectly learn the analytical mapping before any stochastic noise is introduced.}
  \label{fig:fd_predictions}
\end{figure}

The Mean Relative Error (MRE) on this validation set dropped to 1.7\% for absorption ($\mu_a$) and remained well below 5\% for reduced scattering ($\mu_s'$), with success rates (errors $< 10\%$) exceeding 95\%. Figure \ref{fig:fd_predictions} illustrates the predictions of the model against the ground truth for the $\DFD$ test set. The tight clustering along the identity line confirms that the chosen recurrent architecture is perfectly capable of disentangling the non-linear temporal dynamics of $\mu_a$ and $\mu_s'$ when the underlying physical model is smooth and geometrically simplified.

\subsection{Hyperparameter Sensitivity and Network Optimization}
\label{sec:res_hyperparams}

To ensure the robustness of the architecture and avoid suboptimal convergence, an extensive hyperparameter grid search was conducted. We analyzed the influence of the network capacity (number of hidden layers and hidden size), regularization parameters (dropout), the learning rate, and the relative weighting of the $\mu_a$ loss versus the $\mu_s'$ loss. 

\begin{figure}[htbp]
  \centering
  \includegraphics[width=0.49\textwidth]{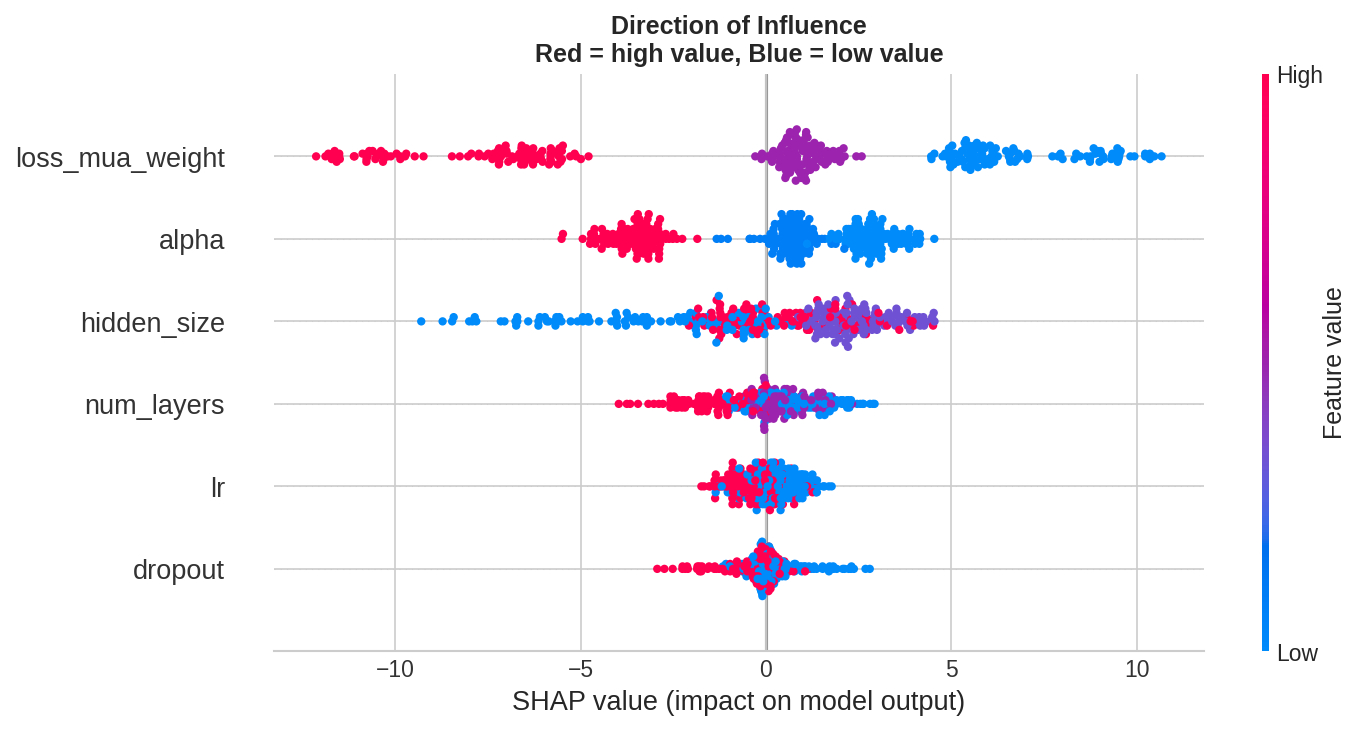}
  \includegraphics[width=0.49\textwidth]{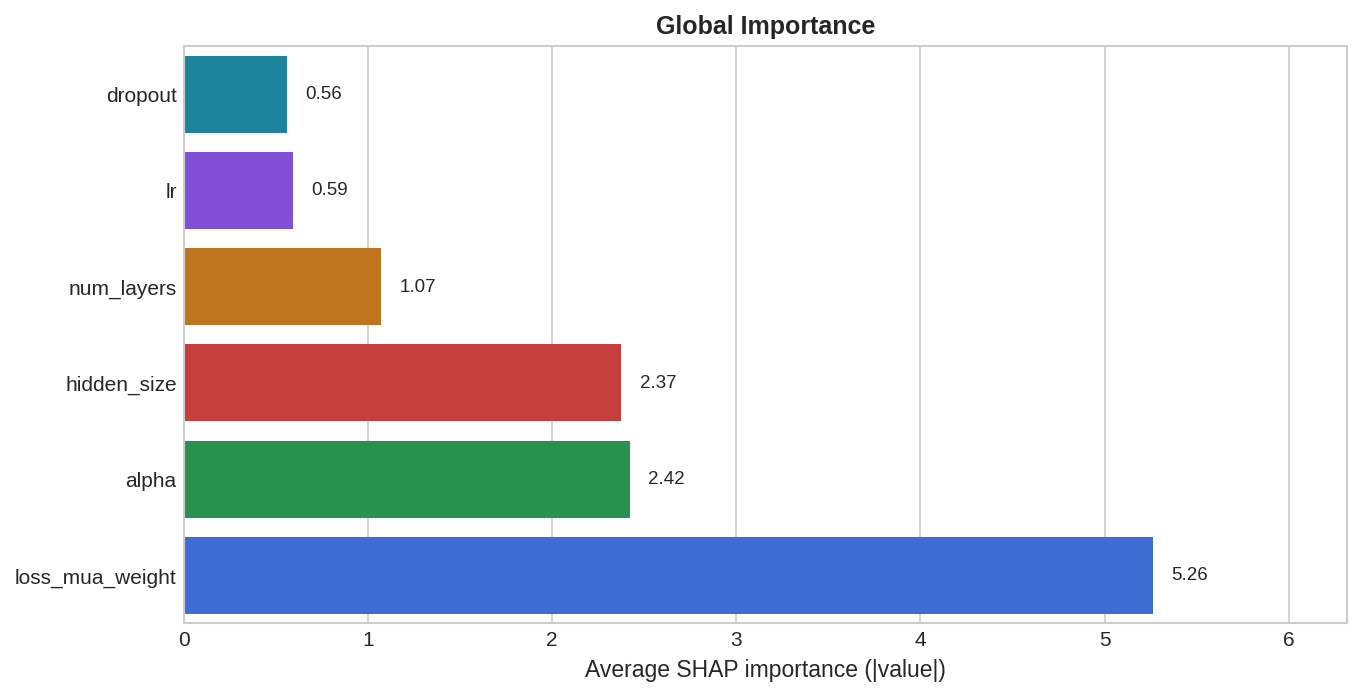}
  \caption{SHAP summary plot illustrating the impact of hyperparameters on the model's success rate. Features are ranked by importance (top to bottom). The $x$-axis shows the marginal impact (SHAP value) on performance, while the color gradient indicates the hyperparameter's actual value (from low in blue, to high in red).}
  \label{fig:hyperparameters}
\end{figure}

To rigorously quantify the influence of each hyperparameter on the model's final performance, a SHapley Additive exPlanations (SHAP) analysis was performed on the grid search results.
SHAP values provide a unified measure of feature importance by calculating the marginal contribution of each parameter to the model's success rate \cite{lundberg2017}. As depicted in Figure \ref{fig:hyperparameters}, the SHAP summary plot ranks the hyperparameters by their overall importance from top to bottom. The horizontal axis represents the SHAP value, which indicates the magnitude and direction of a parameter's impact on the success rate. Furthermore, the color scale denotes the actual tested value of the hyperparameter (e.g., red for high values, blue for low values). 

Interpreting this plot reveals distinct sensitivity profiles. The learning rate emerges as the most critical factor governing convergence stability. The SHAP distribution shows that larger learning rates ($> 10^{-3}$) strongly correlate with negative impacts on the fine-tuning phase, inducing severe oscillations on the noisy $\DMC$ dataset. Reducing the learning rate to $5 \times 10^{-4}$ was imperative to prevent the network from overfitting the stochastic shot noise and is associated with the highest positive SHAP values.

Furthermore, moderately increasing the network capacity (e.g., $\texttt{hidden\_size} = 32$, $\texttt{num\_layers} = 2$) paired with a strict dropout rate of 0.2 provided the optimal balance between feature extraction capability and regularization. Finally, adjusting the loss weight to slightly penalize $\mu_a$ errors more heavily (e.g., a factor of 2.0) consistently yielded positive marginal contributions, helping to synchronize the convergence of both optical parameters despite their different temporal sensitivities.

To further illustrate the network's sensitivity, Table \ref{tab:top5_hyperparams} presents the top 5 hyperparameter (among 324) configurations resulting from the grid search, ranked by their success rate on the reduced scattering coefficient ($\mu_s'$). The optimal configuration confirms the insights derived from the SHAP analysis: a moderate capacity combined with strict dropout and a highly controlled learning rate is essential to successfully handle stochastic noise.

\begin{table}[htbp]
\centering
\caption{Top 5 hyperparameter configurations yielding the highest success rates for $\mu_s'$ extraction during the grid search. The results highlight the preference for moderate network capacity and strict regularization to handle stochastic noise.}
\label{tab:top5_hyperparams}
\begin{tabular}{@{}ccccccc@{}}
\toprule
\textbf{Noise ($\alpha$)} & \textbf{Hidden Size} & \textbf{Layers} & \textbf{Dropout} & \textbf{L.R.} & \textbf{Penal. on $\mu_a$} & \textbf{Success ($\mu_s'$)} \\ \midrule
0.0001 & 32 & 2 & 0.2 & 0.0005 & 2.0 & \textbf{92.4\%} \\
0.0001 & 64 & 1 & 0.2 & 0.0010 & 2.0 & 91.9\% \\
0.0001 & 32 & 2 & 0.1 & 0.0005 & 2.0 & 91.7\% \\
0.0001 & 32 & 3 & 0.1 & 0.0010 & 2.0 & 91.5\% \\
0.0010 & 32 & 1 & 0.2 & 0.0010 & 2.0 & 91.2\% \\ \bottomrule
\end{tabular}
\end{table}

\subsection{Quantifying the Domain Gap}
\label{sec:res_gap}

The core challenge in time-resolved optical inversion is the transition from theoretical models to realistic stochastic measurements. To quantify this analytical-to-stochastic domain gap, we performed a direct inference: the network, trained exclusively on the deterministic $\DFD$ dataset, was tasked with predicting the optical properties of the stochastic 3D Monte Carlo dataset ($\DMC$).

As anticipated by the PCA and temporal correlation analyses (Section \ref{sec:gap_analysis}), this direct application resulted in a catastrophic drop in performance. The network severely misinterpreted the stochastic TPSFs, yielding an MRE of over 60\% for $\mu_a$ and exceeding 150\% for $\mu_s'$. More importantly, the predictions exhibited a massive positive systematic bias (+54.7\% for $\mu_a$ and >200\% for $\mu_s'$). This confirms that the analytical diffusion approximation systematically underestimates the time-of-flight of photons by neglecting 3D boundary losses and the highly directional behavior of early photons, leading the network to artificially inflate the predicted optical properties to compensate for the delayed optical response of the MC data.

\subsection{Bridging the Gap via Transfer Learning}
\label{sec:res_transfer}

Following the physics-informed transfer learning strategy (Section \ref{sec:transfer_learning}), the pre-trained model was fine-tuned using the stochastic target dataset $\DMC$. The dataset was split into training (approx. 2,700 samples), validation (350 samples), and testing (350 samples) subsets to prevent data leakage. 

The fine-tuning phase dramatically reversed the effects of the domain gap. Table \ref{tab:results_summary} summarizes the global performance before and after transfer learning on the Monte Carlo test set.

\begin{table}[htbp]
\centering
\caption{Ablation Study and Global Performance Report on the 3D Monte Carlo Test Set. The table compares direct inference from the analytical prior, training entirely from scratch on the limited MC dataset, and our proposed physics-informed transfer learning. The fine-tuned model significantly outperforms the others, particularly by completely neutralizing the analytical systematic bias and stabilizing predictions.}
\label{tab:results_summary}
\resizebox{\textwidth}{!}{%
\begin{tabular}{@{}l|ccc|ccc@{}}
\toprule
\textbf{Metric} & \multicolumn{3}{c|}{\textbf{Absorption ($\mu_a$)}} & \multicolumn{3}{c}{\textbf{Reduced Scattering ($\mu_s'$)}} \\
 & \textit{Direct (FD)} & \textit{Scratch (MC)} & \textbf{\textit{Fine-Tuned}} & \textit{Direct (FD)} & \textit{Scratch (MC)} & \textbf{\textit{Fine-Tuned}} \\ \midrule
\textbf{MRE} & 64.7\% & 12.8\% & \textbf{10.8\%} & 210.0\% & 22.9\% & \textbf{13.5\%} \\
\textbf{Bias} & +54.7\% & +4.53\% & \textbf{+1.3\%} & +208.3\% & +10.78\% & \textbf{+1.5\%} \\
\textbf{Std. Dev.} & 87.7\% & 20.7\% & \textbf{13.8\%} & 173.8\% & 46.8\% & \textbf{17.3\%} \\ \bottomrule
\end{tabular}%
}
\end{table}

To rigorously validate the necessity of the transfer learning approach, we conducted an ablation study by training an identical Bi-LSTM architecture entirely from scratch using the MC dataset. To ensure a strictly fair comparison, this baseline model was trained on the exact same augmented training subset (5,920 samples) and evaluated on the same isolated test subset (740 samples) as the fine-tuned model. 

As shown in Table \ref{tab:results_summary}, while the from-scratch model manages to partially reduce the massive analytical biases (bringing them down to +4.5\% for $\mu_a$ and +10.8\% for $\mu_s'$), it struggles significantly with overall accuracy and stability. Without the physical prior established by the analytical dataset, the from-scratch network is highly susceptible to stochastic shot noise despite the data augmentation and Savitzky-Golay filtering. It attempts to memorize the random fluctuations of the dataset rather than learning the underlying physics, leading to severe predictive instability. This is evidenced by the massive standard deviation of the relative errors, which spikes to 46.6\% for the reduced scattering coefficient.

The adaptation of our physics-informed fine-tuning overcomes these limitations entirely. It allows the network to achieve a highly competitive MRE of 10.8\% for absorption and 13.5\% for reduced scattering, while dramatically improving the success rates. By initializing the network with a theoretical baseline, the fine-tuning phase acts as a targeted recalibration rather than a blind optimization. It successfully shifts the network's temporal focus to account for 3D boundary losses and the sub-diffusive regime, thereby virtually neutralizing the systematic bias (+1.3\% and +1.5\%) without sacrificing the stability of the predictions (Std. Dev. < 19\%).

\begin{figure}[htbp]
  \centering
  \includegraphics[width=0.48\textwidth]{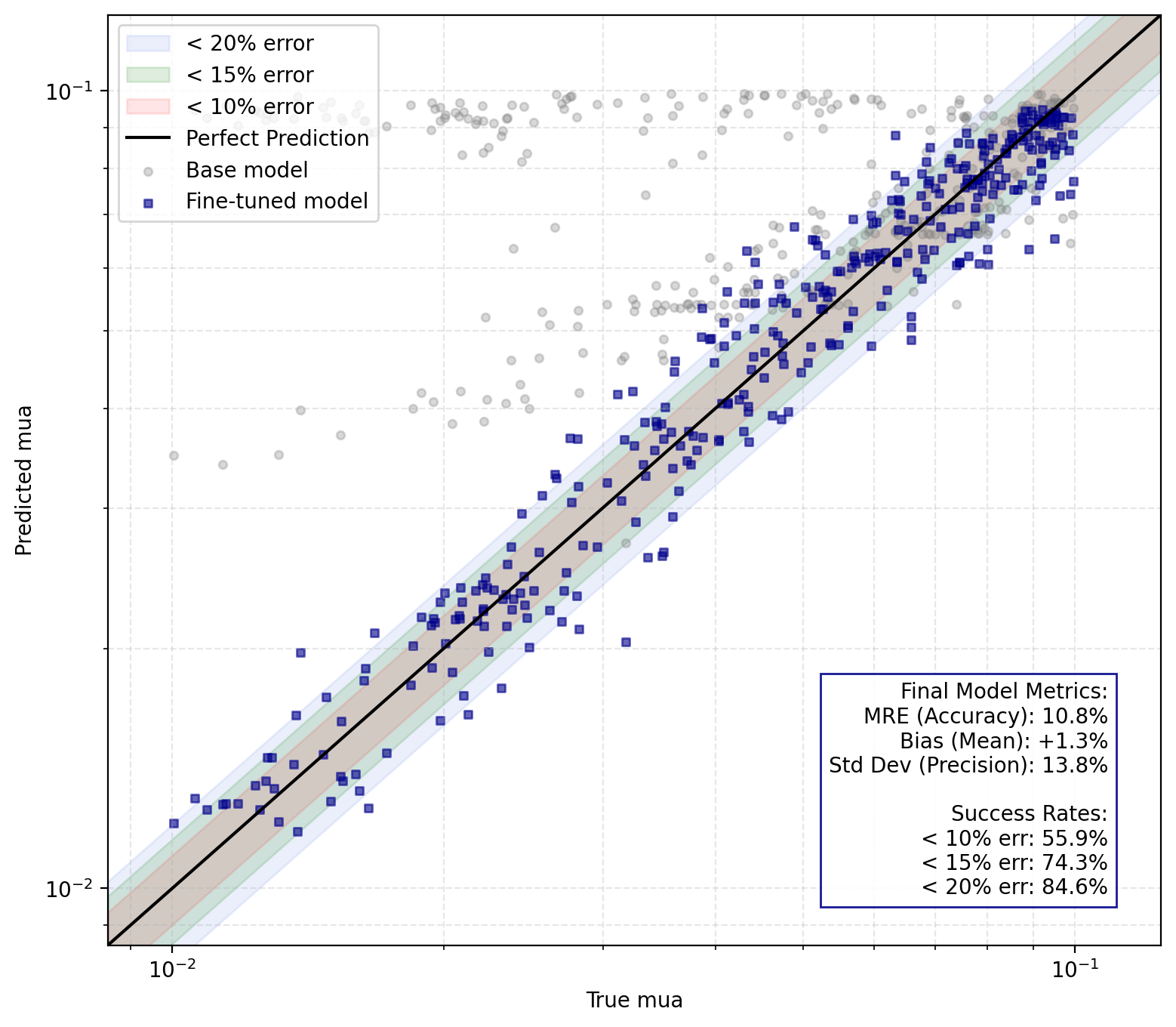}
   \includegraphics[width=0.48\textwidth]{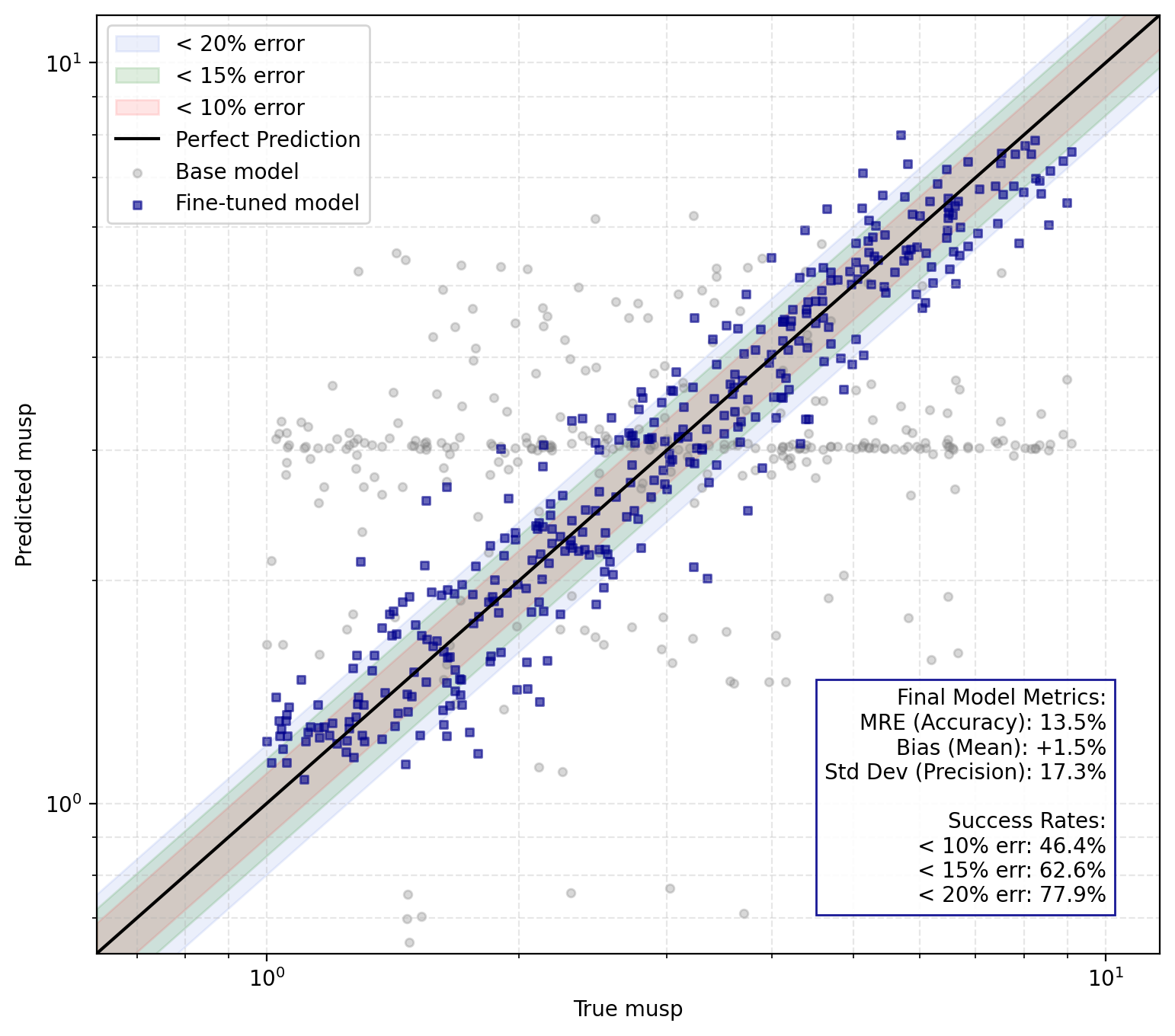}
  \caption{Scatter plots of the predicted optical properties (in blue dots) versus the ground truth Monte Carlo values (in gray dots)  after transfer learning. The solid line represents the ideal prediction ($y = x$). The network demonstrates strong linearity across the parameter space.}
  \label{fig:scatter_predictions}
\end{figure}

\subsection{Discussion}
\label{sec:discussion}

The results validate the premise that large stochastic datasets are not strictly necessary if a robust theoretical prior is established beforehand. By using only 3,700 stochastic samples, our method achieves accuracy levels comparable to recent deep learning studies that required over 100,000 Monte Carlo simulations \cite{li2010, fussen2025}. 

It is noteworthy that the prediction of the reduced scattering coefficient ($\mu_s'$) remains slightly more challenging (higher MRE and Std. Dev.) than the absorption coefficient ($\mu_a$). This behavior aligns with the physics of time-resolved measurements. As demonstrated in Figure \ref{fig:temporal_correlation}, $\mu_a$ governs the long, exponential tail of the TPSF, providing the network with a large temporal window rich in stable features. Conversely, $\mu_s'$ primarily dictates the position and width of the early-time peak (the ballistic and snake photon regimes). This temporal region spans only a few bins and is highly susceptible to stochastic shot noise, making feature extraction inherently more difficult. Despite this physical constraint, the Dual-Head Bi-LSTM successfully mitigates parameter crosstalk and delivers robust, near-instantaneous inversions.

\section{Conclusion and Perspectives}
\label{sec:conclusion}

In this study, we addressed the fundamental trade-off between computational efficiency and physical accuracy in time-resolved optical property inversion. We introduced a highly data-efficient, physics-informed transfer learning framework centered on a Dual-Head Bidirectional LSTM architecture. By reframing the traditional scalar regression as a dual-classification task via a Cross-Entropy loss function, the network not only predicts the absorption ($\mu_a$) and reduced scattering ($\mu_s'$) coefficients but also outputs probability distributions that inherently quantify predictive uncertainty.

Our statistical analyses, including PCA and temporal correlation, rigorously demonstrated the structural domain gap between fast 2D deterministic approximations and realistic 3D stochastic photon transport. We showed that analytical models systematically fail to capture the highly directional sub-diffusive regime and 3D boundary losses, leading to severe predictive biases (+54\% for $\mu_a$ and >200\% for $\mu_s'$). 

By leveraging a two-step transfer learning strategy, the Bi-LSTM was first primed with a strong physical prior using 7,441 inexpensive analytical profiles, and subsequently fine-tuned on a highly restricted set of only 3,700 stochastic Monte Carlo simulations. This synergistic approach successfully bridged the domain gap, virtually eliminating the systematic analytical bias (reduced to under 1.5\%) and achieving a highly competitive Mean Relative Error of 10.8\% for $\mu_a$ and 13.5\% for $\mu_s'$. Crucially, this methodology slashes the exorbitant data requirements—often exceeding 100,000 stochastic samples in recent literature—making high-fidelity deep learning models vastly more accessible. Furthermore, the inference time of the trained network is nearly instantaneous, enabling true real-time characterization of turbid media.

Future work will focus on validating this architecture against experimental measurements from the TROT device \cite{pallares2021}.

\section{Reproducibility}

The source code for the training, Monte Carlo 3D simulation, and inversion components is publicly accessible at \url{https://gitlab.math.unistra.fr/aghili/anatrot_mc}. The repository also includes the three models referenced in Table \ref{tab:results_summary}, provided in PyTorch \texttt{pth} format, along with the two datasets, $\DFD$ and $\DMC$, in Numpy binary format \texttt{npy}. By following the instructions outlined in the repository, users can easily reproduce all figures presented in the study.

\section*{Declaration of generative AI in the manuscript preparation process.}

During the preparation of this work the author(s) used Google Gemini 3.1 (mostly) and Claude 4.6 in order to implement ideas into Python code, generate plot various scripts, translate and enhance texts from French to English. After using this tool/service, the author(s) reviewed and edited the content as needed and take(s) full responsibility for the content of the published article.

\section*{Acknowledgments}

The authors wish to express their gratitude to the Interdisciplinary Thematic Institute IRMIA++, which operates under the ITI 2021-2028 program of the University of Strasbourg, CNRS, and Inserm, receiving support from IdEx Unistra (ANR-10-IDEX-0002) and the SFRI-STRAT'US project (ANR-20-SFRI-0012) within the framework of the French Investments for the Future Program. Additionally, the authors extend their appreciation to E. Franck for his insightful discussions. Finally, they acknowledge the High Performance Computing Center of the University of Strasbourg for its support in providing scientific assistance and access to computing resources, with part of these resources funded by the Equipex Equip@Meso project (\textit{Programme Investissements d'Avenir}) and the CPER Alsacalcul/Big Data.

\bibliographystyle{plain}
\bibliography{paper}

\end{document}